\input amstex

\documentstyle{amsppt}

\def\normo#1{\left\|#1\right\|}
\def\modo#1{\left|#1\right|}

\def\Z{\Bbb Z}
\def\R{\Bbb R}
\def\T{\Bbb T}

\def\N{\Bbb N}

\def\Re{\mathop{\hbox{\rm Re}}}

\def\lub{\mathop{\hbox{\rm l.u.b.}}}

\topmatter

\title
Stability and Dichotomy of Positive Semigroups on $L_p$
\endtitle
\author
Stephen Montgomery-Smith
\endauthor
\thanks
Research supported
in part by N.S.F.\ Grant D.M.S.\ 9201357.
\endthanks
\address
Math.\ Dept.,
University of Missouri, Columbia MO 65211, U.S.A.
\endaddress


\abstract
A new proof of a result of Lutz Weis is given, that states that the
stability of a positive 
strongly continuous semigroup $(e^{tA})_{t \ge 0}$\ on $L_p$\ 
may be determined by the quantity $s(A)$.  We also give an example
to show that the dichotomy of the semigroup may not always be
determined by the spectrum $\sigma(A)$.
\endabstract

\endtopmatter

Consider a strongly continuous semigroup $(e^{tA})_{t\ge0}$\ acting on a 
Banach
space $X$\ with unbounded generator $A$.  It has long been known that
the spectral mapping theorem
$e^{t \sigma(A)} = \sigma(e^{tA})\setminus \{0\}$\ 
does not necessarily hold.  (Here
$\sigma(A)$\ denotes the spectrum of an operator $A$.)  Indeed, let
$s(A) = \sup\Re(\sigma(A))$, and let
$\omega(A) 
= \sup\Re(\log(\sigma(e^A))) 
= \inf\{\lambda : \normo{e^{tA}} \le M_\lambda e^{\lambda t}\}$.
Then there are examples of semigroups for which $s(A) \ne \omega(A)$\
(see \cite{N}).  

The purpose of this paper is to give one situation in which it is 
true that $s(A) = \omega(A)$.  This next result has already been proved
by Lutz Weis \cite{We}.  
We will give a different, shorter proof.
We refer the reader
to \cite{We} for a history of the problem.

\proclaim{Theorem 1} Let $e^{tA}$\ be a strongly continuous positive
semigroup on $L_p(\Omega,{\Cal F},\mu)$, where $(\Omega,{\Cal F},\mu)$\
is a sigma-finite measure space, and $1 \le p < \infty$.
Then $\omega(A) = s(A)$.
\endproclaim

In order to show this result, we will make use of the following lemmas.
The first result may be derived from \cite{C}, Theorem~7.4 
(the reader may like to know
that a proof of the `Pringsheim-Landau Theorem' used in \cite{C}
may be found on page~59
of \cite{Wi}).

\proclaim{Lemma 2}  Let $e^{tA}$\ be a strongly continuous positive
semigroup on a Banach
lattice $X$, and let $g \in X$.  Then for any $\lambda > s(A)$\
we have that
$$ (\lambda - A)^{-1} g = \int_0^\infty e^{s(A-\lambda)} g \, ds .$$
Here the right hand side is taken in the sense of an improper integral.
\endproclaim

The next result may be found in \cite{LM1} and \cite{LM2}.

\proclaim{Lemma 3}  Let $e^{tA}$\ be a strongly continuous semigroup on
a Banach space $X$, and let $1 \le p < \infty$.
Then $1 \notin \sigma(e^{2\pi A})$\ 
if and only if $i \Z \cap \sigma(A) = \emptyset$\
and there is a constant $c > 0$\ such that for any $v_{-n}$, 
$v_{-n+1},\dots$,\ $v_n \in X$\ we have
$$ \int_0^{2\pi} \normo{\sum_{k=-n}^n (ik - A)^{-1} v_k e^{ikt}}^p \, dt
   \le
   c^p
   \int_0^{2\pi} \normo{\sum_{k=-n}^n v_k e^{ikt}}^p \, dt .$$
\endproclaim

For the next result, we specialize to a Banach lattice of functions
on a sigma-finite measure space.
In fact, this is really no loss of generality, and
the interested reader
should find no trouble making sense of this result for a general
Banach lattice by applying the ideas in \cite{LT} Chapter~1.4.

\proclaim{Lemma 4}  Let $P$\ be a positive operator on $X$, a 
Banach lattice of functions on a sigma-finite measure space,
such that $\modo{g} \le f \in X$\ implies that $g \in X$.
Let $1 \le p < \infty$.  
If $f:[0,2\pi] \to X$\ is a measurable, simple function,
then
$$ \left( \int_0^{2\pi} \modo{P(f(t))}^p \, dt \right)^{1/p}
   \le
   P\left(\left(\int_0^{2\pi} \modo{f(t)}^p \, dt \right)^{1/p} \right) .$$
\endproclaim

\demo{Proof}
Let us set $f = \sum_{k=1}^n v_k \chi_{A_k}$, 
where $v_k \in X$, and the sets 
$A_k \subseteq [0,2\pi]$\
are disjoint.
Then, letting $f_k = v_k \modo{A_k}^{1/p}$, the result 
reduces to showing that
$$ \left( \sum_{k=1}^n \modo{P(f_k)}^p \right)^{1/p} 
   \le
   P\left(\left( \sum_{k=1}^n \modo{f_k}^p \right)^{1/p} \right).$$
However, we know that
$$ \left( \sum_{k=1}^n \modo{f_k}^p \right)^{1/p}
   =
   \lub_{\sum \modo{a_k}^q \le 1} \sum_{k=1}^n
   \Re(a_k f_k) .$$
Here, $\lub$\ denotes the least upper bound in the lattice.
Now, since $P$\ is positive, we have that
$$ P\left(\lub_{\sum \modo{a_k}^q \le 1} \sum_{k=1}^n
   \Re(a_k f_k) \right) $$ 
is an upper bound for $\sum_{k=1}^n
\Re(a_k P(f_k))$\ whenever $\sum \modo{a_k}^q \le 1$.
Hence
$$ \eqalignno{
   \left( \sum_{k=1}^n \modo{P(f_k)}^p \right)^{1/p}
   &=
   \lub_{\sum \modo{a_k}^q \le 1} \sum_{k=1}^n
   \Re(a_k P(f_k)) \cr
   &\le
   P\left(
   \lub_{\sum \modo{a_k}^q \le 1} \sum_{k=1}^n
   \Re(a_k f_k) \right) \cr
   &=
   P\left(\left( \sum_{k=1}^n \modo{f_k}^p \right)^{1/p} \right) .\cr}$$
\enddemo

\demo{Proof of Theorem~1}
It is well known that $s(A) \le \omega(A)$ (see \cite{N}).  
Thus
by simple rescaling arguments, we see that
it is sufficient to show that if $s(A) < 0$, then 
$\T \cap \sigma(e^{2\pi A}) = \emptyset$.

We will show, under the assumption that $s(A)<0$, that if 
$f:\R \to L_p$\ is a bounded, measurable
function that is periodic with period $2\pi$, then for each $N>0$\
we have
$$ \left( \int_0^{2\pi} \normo{\int_0^N e^{sA} f(t-s) \, ds }_{L_p}^p
   \, dt \right)^{1/p} 
   \le \normo{A^{-1}} \left(\int_0^{2\pi} \normo {f(t)}_{L_p}^p
      \, dt \right)^{1/p} . $$
In order to show this, we may assume without loss of generality that
$f$\ restricted to $[0,2\pi]$\ is a simple function.  Fix $N > 0$.
By the positivity of $e^{sA}$, and Fubini's Theorem, we have that
$$ \eqalignno{
   \left( \int_0^{2\pi} \normo{\int_0^N e^{sA} f(t-s) \, ds }_{L_p}^p
   \, dt \right)^{1/p} 
   &\le
   \left( \int_0^{2\pi} \normo{\int_0^N e^{sA} \modo{f(t-s)} \, ds }_{L_p}^p
   \, dt \right)^{1/p} \cr
   &=
   \normo{ \left( \int_0^{2\pi} \left(\int_0^N e^{sA} \modo{f(t-s)} 
   \, ds \right)^p\, dt \right)^{1/p}}_{L_p} . \cr}$$
By the integral version of 
Minkowski's Theorem (see \cite{HLP}, Section 203), it follows that
for each $\omega \in \Omega$
$$ \eqalignno{
   \left( \int_0^{2\pi} \left(\int_0^N e^{sA} \modo{f(t-s)(\omega)}
   \,ds \right)^p 
   \, dt \right)^{1/p}
   &\le
   \int_0^N \left( \int_0^{2\pi} \left( e^{sA} \modo{f(t-s)(\omega)} \right)^p 
   \, dt \right)^{1/p} \,ds \cr
   &=
   \int_0^N \left( \int_0^{2\pi} \left( e^{sA} \modo{f(t)(\omega)} \right)^p 
   \, dt \right)^{1/p} \,ds .\cr }$$
Finally, from Lemma~4, we see that
$$ \left( \int_0^{2\pi} (e^{sA} \modo{f(t)})^p 
   \, dt \right)^{1/p}
   \le
   e^{sA} \left( \int_0^{2\pi} \modo{f(t)}^p 
   \, dt \right)^{1/p} .$$
Putting all of these together, and applying Lemma~2, we obtain
$$ \eqalignno{
   \left( \int_0^{2\pi} \normo{\int_0^N e^{sA} f(t-s) \, ds }_{L_p}^p
   \, dt \right)^{1/p} 
   &\le \normo{ \int_0^N e^{sA} \left( \int_0^{2\pi} \modo{f(t)}^p \, dt
   \right)^{1/p} \, ds }_{L_p} \cr
   &\le \normo{ \int_0^\infty e^{sA} \left( \int_0^{2\pi} \modo{f(t)}^p \, dt
   \right)^{1/p} \, ds }_{L_p} \cr
   &= \normo{ A^{-1} \left( \int_0^{2\pi} \modo{f(t)}^p \, dt \right)^{1/p}
   }_{L_p} \cr
   &\le \normo{A^{-1}} \normo{ \left( 
   \int_0^{2\pi} \modo{f(t)}^p \, dt \right)^{1/p}
   }_{L_p} \cr
   &= \normo{A^{-1}} \left(\int_0^{2\pi} \normo{f(t)}_{L_p}^p
      \, dt \right)^{1/p} , \cr } $$
where the last equality uses Fubini's theorem.
Now, if $f(t) = e^{i \beta t} \sum_{k=-n}^n v_k e^{ikt}$\ for 
some $\beta \in \R$, 
then by Lemma~2, we see that
$$ \int_0^N e^{sA} f(t-s) \, ds \to 
   \sum_{k=-n}^n (ik + i\beta - A)^{-1} v_k e^{ikt} $$
uniformly in $t$\ as $N\to\infty$.
Hence by Lemma~3 it follows that $e^{i\beta} \notin \sigma(e^{2\pi A})$.
\enddemo

One might
conjecture that the spectrum of the generator of a positive semigroup 
$e^{tD}$\ on an $L_p$\ space might
characterize the dichotomy of the semigroup, that is, if $a$\ is any real 
number, then
$(a + i\R) \cap \sigma(D) = \emptyset$\ if and only if
$e^{ta} \T \cap \sigma(e^{tD}) = \emptyset$.
However, this is not the case, as the next result shows.

\proclaim{Theorem 5}  There is a positive semigroup $e^{tD}$\ acting on 
an $L_2$\ space such that $(1 + i\R) \cap \sigma(D) = \emptyset$,
but $e^{2\pi} \in \sigma(e^{2\pi D})$.
\endproclaim

\demo{Proof}
For each $M \in \N$,
let $C_M$\ be the contraction acting on $\ell_2^M$\ by the matrix
$$ C_M = \left[ \matrix
   0 & 1 & 0 & 0 & \cdots & 0 \cr
   0 & 0 & 1 & 0 & \cdots & 0 \cr
   0 & 0 & 0 & 1 & \cdots & 0 \cr
   0 & 0 & 0 & 0 & \cdots & 0 \cr
   \vdots & \vdots & \vdots & \vdots & & \vdots \cr
   0 & 0 & 0 & 0 & \cdots & 0 \cr 
   \endmatrix \right] .$$

Note that if $\lambda \ne 0$, then
$$ (\lambda - C_M)^{-1} 
= \sum_{j=0}^{M-1} \lambda^{-1-j} C_M^{j} 
= \left[ \matrix
\lambda^{-1} & \lambda^{-2} & \lambda^{-3} & \lambda^{-4} 
             & \cdots & \lambda^{-M}   \cr
0            & \lambda^{-1} & \lambda^{-2} & \lambda^{-3} 
             & \cdots & \lambda^{-M+1} \cr
0            & 0            & \lambda^{-1} & \lambda^{-2} 
             & \cdots & \lambda^{-M+2} \cr
0            & 0            & 0            & \lambda^{-1} 
             & \cdots & \lambda^{-M+3} \cr
\vdots       & \vdots       & \vdots       & \vdots       
             &        & \vdots         \cr
0            & 0            & 0            & 0            
             & \cdots & \lambda^{-1}   \cr
\endmatrix \right]  .$$
Thus, if $\modo \lambda = 1$, then $\normo{(\lambda - C_M)^{-1}}
\ge \sqrt M$.  Also, if $\modo{\lambda} > 1$, then
$\normo{(\lambda - C_M)^{-1}}
\le \sum_{j=0}^{M-1} \modo\lambda^{-1-j} \le 1/(\modo\lambda - 1)$.
In particular, if $\modo\lambda \ge 2$, then 
$\normo{(\lambda - C_M)^{-1}} \le 1$.

Also note that
$$ e^{tC_M} = \left[ \matrix
1 & t & t^2/2 & t^3/6 & \cdots & t^{M-1}/(M-1)! \cr
0 & 1 & t     & t^2/2 & \cdots & t^{M-2}/(M-2)! \cr
0 & 0 & 1     & t     & \cdots & t^{M-3}/(M-3)! \cr
0 & 0 & 0     & 1     & \cdots & t^{M-4}/(M-4)! \cr
\vdots & \vdots & \vdots & \vdots & & \vdots    \cr
0 & 0 & 0     & 0     & \cdots & 1              \cr
\endmatrix \right] .$$
Thus we see that $e^{tC_M}$\ is a positive operator.
Clearly $\normo{e^{tC_M}} \le e^{t\normo{C_M}}
\le e^t$.

Consider the positive semigroup acting on $L_2([0,2\pi])$\ by
$$ e^{t A_M}f(x) = (e^{4t}-1) \int_0^{2\pi} f(x) \, {dx\over 2\pi}
   + f(x+Mt) ,$$
so that its generator is the closure of
$$ A_M f(x) = 4 \int_0^{2\pi} f(x) \, {dx\over 2\pi}
   + M {d \over dx} f(x) .$$
Note that $\normo{e^{t A_M}} \le e^{4t}$.

Now consider the positive
semigroup $e^{t B_M} = e^{tA_M} \otimes e^{tC_M}$\ 
acting on
$$ X_M = L_2([0,2\pi]) \otimes \ell_2^M = 
   L_2([0,2\pi] \times \{1,2,\dots,M\}) ,$$
We see that this semigroup 
is generated by $B_M = A_M \otimes I + I \otimes C_M$.
Also, $\normo{e^{tB_M}} \le e^{5t}$.

Consider a typical element of $X_M$\ given by 
$f(x) = \sum_{n=-\infty}^\infty v_n e^{i n x} \in X_M$, where
$v_n \in \ell_2^M$, and $\normo f_{X_M}^2 = 2\pi \sum_{n=-\infty}^\infty
\normo{v_n}_2^2$.
If $\lambda \ne 4$\ and $\lambda \notin M\Z\setminus\{0\}$, then $\lambda 
\notin \sigma(B_M)$,
and
$$ (\lambda - B_M)^{-1} f(x) =
   (\lambda - 4 - C_M)^{-1} v_0
   + \sum_{n\ne0} (\lambda - i n M - C_M)^{-1} v_n e^{i n x} .$$
Thus
$$ \normo{(\lambda - B_M)^{-1}} =
   \max\left\{ \normo{(\lambda - 4 - C_M)^{-1}} ,
   \sup_{n\ne0} \normo{(\lambda - i n M - C_M)^{-1}} \right\} .$$
In particular, if $\Re(\lambda) = 1$\ and $\modo{\lambda} \le M - 2$,
then $\normo{(\lambda - B_M)^{-1}} \le 1$, whereas if 
$\lambda = 1 + i M$,
then $\normo{(\lambda - B_M)^{-1}} \ge \sqrt M$.

Now consider the semigroup $e^{tD} =
\bigoplus\limits_{M=1}^\infty e^{tB_M}$\ acting on
$$ \bigoplus_{M=1}^\infty X_M
   = L_2\left(\bigvee_{M=1}^\infty 
   \bigl([0,2\pi] \times \{1,2,\dots,M\}\bigr)\right) .$$
Note that $e^{tD}$\ really is a strongly continuous semigroup,
with $\normo{e^{tD}} \le e^{5t}$.
The generator $D$\ is the closure of $\bigoplus\limits_{M=1}^\infty B_M$, and
hence its resolvent set consists of those $\lambda$\ such that
$$ \normo{(\lambda - D)^{-1}} 
   = \sup_{M \ge 1} \normo{(\lambda - B_M)^{-1}}  < \infty,$$
that is, $\sigma(D) \subseteq \{z:\modo{z-4} \le 1\} \cup i\Z \setminus\{0\}$.
In particular, if $\Re(\lambda) = 1$, then $\lambda \notin \sigma(D)$.  
However,
$ \sup_{\lambda\in1+i\Z} \normo{(\lambda - D)^{-1}} = \infty$, and
hence, by Gerhard's Theorem (see \cite{N}, p.~95), $e^{2\pi} 
\in \sigma(e^{2\pi D})$.
\enddemo

\refindentwd = 10 true mm

\enddocument